\documentclass[12pt]{smfart}

\usepackage{smfthm}
\usepackage{amssymb}
\usepackage[francais]{babel}
\usepackage[applemac]{inputenc}

\newcommand{\Qp}{\mathbf{Q}_p}
\newcommand{\Fp}{\mathbf{F}_p}
\newcommand{\Zp}{\mathbf{Z}_p}
\newcommand{\ZZ}{\mathbf{Z}}
\newcommand{\QQ}{\mathbf{Q}}
\newcommand{\FF}{\mathbf{F}}
\newcommand{\RR}{\mathbf{R}}
\newcommand{\OO}{\mathcal{O}}
\newcommand{\MM}{\mathfrak{m}}
\newcommand{\Qpbar}{\overline{\mathbf{Q}}_p}
\newcommand{\Qbar}{\overline{\mathbf{Q}}}
\renewcommand{\phi}{\varphi}

\renewcommand{\geq}{\geqslant}
\renewcommand{\leq}{\leqslant} 
\newcommand{\Gal}{\mathrm{Gal}}
\newcommand{\Fil}{\mathrm{Fil}}
\newcommand{\bst}{\mathbf{B}_{\mathrm{st}}} 
\newcommand{\bcris}{\mathbf{B}_{\mathrm{cris}}} 
\newcommand{\calE}{\mathcal{E}}
\newcommand{\calR}{\mathcal{R}}
\newcommand{\calS}{\mathcal{S}}
\newcommand{\calC}{\mathcal{C}}
\newcommand{\bdr}{\mathbf{B}_{\mathrm{dR}}}  
\newcommand{\bfont}{\mathbf{B}}
\newcommand{\dst}{\mathrm{D}_{\mathrm{st}}}
\newcommand{\dpst}{\mathrm{D}_{\mathrm{pst}}}
\newcommand{\dcris}{\mathrm{D}_{\mathrm{cris}}}
\newcommand{\ddr}{\mathrm{D}_{\mathrm{dR}}}
\newcommand{\dfont}{\mathrm{D}}
\newcommand{\smat}[1]{\left( \begin{smallmatrix} #1 \end{smallmatrix} \right)}
\newcommand{\pmat}[1]{\begin{pmatrix} #1 \end{pmatrix}}
\newcommand{\B}{\mathrm{B}}
\newcommand{\K}{\mathrm{K}}
\newcommand{\vp}{\mathrm{val}_p}
\newcommand{\Linv}{\mathcal{L}}
\newcommand{\ab}{\mathrm{ab}}

\newcommand{\Res}{\mathrm{Res}}
\newcommand{\End}{\mathrm{End}}
\newcommand{\Aut}{\mathrm{Aut}}
\newcommand{\Sym}{\mathrm{Sym}}
\newcommand{\GL}{\mathrm{GL}}
\newcommand{\LL}{\mathrm{Lisse}}
\newcommand{\AL}{\mathrm{Alg}}
\newcommand{\Sp}{\mathrm{Sp}}
\newcommand{\St}{\mathrm{St}}
\newcommand{\ind}{\mathrm{ind}}
\newcommand{\an}{\mathrm{an}}
\newcommand{\alg}{\mathrm{alg}}

\newcommand{\irr}{\mathrm{irr}}
\newcommand{\nr}{\mathrm{nr}}
\newcommand{\cris}{\mathrm{cris}}
\newcommand{\st}{\mathrm{st}}
\newcommand{\ssf}{\mathrm{ss}}
\newcommand{\bor}{\mathrm{bor}}
\newcommand{\ngeo}{\mathrm{ng}}
\newcommand{\lisse}{\mathrm{lisse}}
\newcommand{\Mat}{\mathrm{Mat}}
\newcommand{\Ext}{\mathrm{Ext}}
\newcommand{\drig}{\dfont_{\mathrm{rig}}}
\newcommand{\dsharp}{\dfont^\sharp}

\newcommand{\Proj}{\mathbf{P}}
\newcommand{\dbf}{\mathbf{D}}
\newcommand{\Pibf}{\mathbf{\Pi}}
\newcommand{\dcroc}[1]{[\![ #1 ]\!]}

\date{Avril 2010}
\title[La correspondance de Langlands locale $p$-adique pour $\mathrm{GL}_2(\mathbf{Q}_p)$]{La correspondance de Langlands locale $p$-adique pour $\mathrm{GL}_2(\mathbf{Q}_p)$ \\ E\lowercase{xposé} \no 1017 \lowercase{du} S\lowercase{éminaire} B\lowercase{ourbaki}}

\alttitle{The $p$-adic local Langlands correspondence for $\mathrm{GL}_2(\mathbf{Q}_p)$}

\author{Laurent Berger}
\address{UMPA ENS Lyon \\
UMR 5669 du CNRS \\
Universit\'e de Lyon}
\email{laurent.berger@umpa.ens-lyon.fr}
\urladdr{www.umpa.ens-lyon.fr/\~{}lberger/}



\NumberTheoremsIn{subsection}

\begin{document}

\begin{abstract}
La correspondance de Langlands locale $p$-adique pour $\mathrm{GL}_2(\mathbf{Q}_p)$ est une bijection entre certaines repr\'esentations de dimension $2$ de $\mathrm{Gal}(\overline{\mathbf{Q}}_p / \mathbf{Q}_p)$ et certaines repr\'esen\-tations de $\mathrm{GL}_2(\mathbf{Q}_p)$. Cette bijection peut en fait \^etre construite en utilisant la th\'eorie des $(\varphi,\Gamma)$-modules et des r\'esultats d'analyse $p$-adique. On d\'eduit alors des propri\'et\'es de cette construction quelques applications int\'eressantes en arithm\'etique.
\end{abstract}

\begin{altabstract}
The $p$-adic local Langlands correspondence for $\mathrm{GL}_2(\mathbf{Q}_p)$ is a bijection between some $2$-dimensional representations of $\mathrm{Gal}(\overline{\mathbf{Q}}_p / \mathbf{Q}_p)$ and some representations of $\mathrm{GL}_2(\mathbf{Q}_p)$. This bijection can in fact be constructed using the theory of $(\varphi,\Gamma)$-modules and some results of $p$-adic analysis. One then deduces from the properties of this construction some interesting arithmetical applications. 
\end{altabstract}

\maketitle

\setcounter{tocdepth}{2}
\tableofcontents

\setlength{\baselineskip}{18pt}

\section*{Introduction}

La correspondance de Langlands locale $p$-adique pour $\GL_2(\Qp)$ est une bijection entre certaines repr\'esentations de dimension $2$ de $\Gal(\Qpbar/\Qp)$ et certaines repr\'esentations de $\GL_2(\Qp)$. Ces représentations sont à coefficients soit dans une extension finie de $\Fp$ (correspondance en caractéristique $p$), soit dans une extension finie de $\Qp$ (correspondance $p$-adique). 

Dans le cas de la correspondance en caractéristique $p$, on peut faire une liste des objets du côté $\Gal(\Qpbar/\Qp)$ ainsi que  du côté $\GL_2(\Qp)$ et cela permet de définir une bijection numérique dont la construction est donnée dans le \S\ref{repcarp}. Dans le cas de la correspondance $p$-adique, on commence par expliquer comment faire le tri dans les représentations de $\Gal(\Qpbar/\Qp)$ (théorie de Fontaine) et dans celles de $\GL_2(\Qp)$ (représentations {\og admissibles\fg} au sens de Schneider et Teitelbaum). Ensuite, on donne les premiers exemples de correspondance construits par Breuil, ce qui est l'objet du \S\ref{repnul}. C'est en étudiant ces exemples que Colmez a compris comment construire de manière fonctorielle cette correspondance, grâce à la théorie des $(\varphi,\Gamma)$-modules. Cette construction est donnée dans le \S\ref{constcorr} pour les représentations {\og triangulines\fg}. Le \S\ref{constgl} contient la construction générale, ainsi que quelques propriétés de cette correspondance
\begin{enumerate}
\item la compatibilité à la réduction modulo $p$,
\item le lien avec la correspondance de Langlands locale {\og classique\fg}.
\end{enumerate}
Dans le \S\ref{apps}, nous donnons quelques applications de la correspondance, dont
\begin{enumerate}
\item le calcul de la réduction modulo $p$ des représentations cristallines,
\item la démonstration de (nombreux cas de) la conjecture de Fontaine-Mazur.
\end{enumerate}
Ce texte ne dit pas grand chose sur la motivation de ces constructions. Outre la compatibilité avec la correspondance de Langlands locale classique, une propriété très importante de la correspondance $p$-adique est sa réalisation dans la cohomologie complétée des tours de courbes modulaires (travaux en cours de rédaction). Enfin, l'extension de ces constructions à d'autres groupes que $\GL_2(\Qp)$ est particulièrement délicate et fait l'objet de nombreux travaux en cours dont il serait prématuré de parler (voir \cite{BICM}).

\subsection*{Notations}

Dans tout ce texte, $E$ désigne une extension finie de $\Qp$ dont on note $\OO_E$ l'anneau des entiers, $\MM_E$ l'idéal maximal de $\OO_E$ et $k_E$ son corps résiduel. Les corps $E$ et $k_E$ sont les corps des coefficients des représentations que l'on considère.

La théorie du corps de classes local fournit une application $\Qp^\times \to \Gal(\Qpbar/\Qp)^{\ab}$ dont l'image est dense et que l'on normalise en décidant que l'image de $p$ est le frobenius géométrique. Cette application nous permet de considérer les caractères de $\Qp^\times$ à valeurs dans $E^\times$ ou $k_E^\times$ comme des caractères de $\Gal(\Qpbar/\Qp)^{\ab}$. On note $\mu_\lambda$ le caractère de $\Qp^\times$ qui est non-ramifié (c'est-à-dire trivial sur $\Zp^\times$) et qui envoie $p$ sur $\lambda$ et on note $|{\cdot}|$ le caractère $x \mapsto p^{-\vp(x)}$ où $\vp(p)=1$.

\section{Représentations en caractéristique $p$}
\label{repcarp}

Dans cette section, nous donnons la classification des représentations de $\Gal(\Qpbar/\Qp)$ et de $\GL_2(\Qp)$ en caractéristique $p$, afin de définir la correspondance dans ce cas.

\subsection{Représentations de $\Gal(\Qpbar/\Qp)$}
\label{galmodp}

Soit $\Qp^{\nr}$ l'extension maximale non-ramifiée de $\Qp$ de telle sorte que $\Gal(\Qpbar/\Qp^{\nr})$ est le sous-groupe d'inertie $\mathcal{I}_{\Qp}$ de $\Gal(\Qpbar/\Qp)$. Si $n \geq 1$ et $d=p^n-1$, alors $\Qp^{\nr}(p^{1/d})$ est une extension modérément ramifiée de $\Qp^{\nr}$ et l'application $g \mapsto g(p^{1/d})/(p^{1/d})$ définit par réduction modulo $p$ un caractère $\omega_n : \mathcal{I}_{\Qp} \to \FF_{p^n}^\times$ (c'est le caractère {\og de niveau $n$\fg} $\theta_{p^n-1}$ du \S 1.7 de \cite{S72}). Par exemple, $\omega_1$ est la restriction à $\mathcal{I}_{\Qp}$ de $\omega$, la réduction modulo $p$ du caractère cyclotomique.

Si $h \in \ZZ$, alors il existe une unique représentation semi-simple notée $\ind(\omega_n^h)$, de déterminant $\omega^h$ et de restriction à $\mathcal{I}_{\Qp}$ isomorphe à $\omega_n^h \oplus \omega_n^{ph} \oplus \cdots \oplus \omega_n^{p^{n-1}h}$. La représentation $\ind(\omega_n^h)$ est alors $\Fp$-linéaire. Si $\chi : \Gal(\Qpbar/\Qp) \to k_E^\times$ est un caractère, alors on note $\rho(r,\chi)$ la représentation $\ind(\omega_2^{r+1}) \otimes \chi$ qui est absolument irréductible si $r \in \{0,\hdots,p-1\}$. 

\begin{theo}\label{repskegal}
Toute représentation $k_E$-linéaire absolument irréductible de 
\mbox{dimension}~$2$ de $\Gal(\Qpbar/\Qp)$ est isomorphe à $\rho(r,\chi)$ pour un $r \in \{0,\hdots,p-1\}$.
\end{theo}

Toute représentation $k_E$-linéaire semi-simple de dimension $2$ de $\Gal(\Qpbar/\Qp)$ est donc isomorphe (après extension éventuelle des scalaires) à $\rho(r,\chi)$ ou bien à $\omega^r \mu_\lambda \oplus \omega^s \mu_\nu$.

\subsection{Représentations de $\GL_2(\Qp)$}
\label{gl2modp}
Nous donnons à présent la classification des représentations $k_E$-linéaires \emph{lisses} (c'est-à-dire localement constantes) et absolument irréductibles de $\GL_2(\Qp)$ qui admettent un caractère central. On identifie le centre de $\GL_2(\Qp)$ à $\Qp^\times$. Si $r \geq 0$, alors $\Sym^r k_E^2$ est une représentation de $\GL_2(\Fp)$ qui fournit par inflation une représentation de $\GL_2(\Zp)$, et on l'étend à $\GL_2(\Zp)\Qp^\times$ en faisant agir $p$ par l'identité. La représentation 
\[ \ind_{\GL_2(\Zp)\Qp^\times}^{\GL_2(\Qp)} \Sym^r k_E^2 \] 
est l'ensemble des fonctions $f : \GL_2(\Qp) \to \Sym^r k_E^2$ qui sont localement constantes, à support compact modulo $\GL_2(\Zp)\Qp^\times$ et telles que $f(kg)=\Sym^r(k) f(g)$ si 
\mbox{$k \in \GL_2(\Zp)\Qp^\times$} et $g \in \GL_2(\Qp)$. Cet espace est muni de l'action de $\GL_2(\Qp)$ donnée par $(gf)(h)=f(hg)$. L'algèbre de Hecke 
\[ \End_{k_E[\GL_2(\Qp)]} \left( \ind_{\GL_2(\Zp)\Qp^\times}^{\GL_2(\Qp)} \Sym^r k_E^2 \right) \] 
se calcule à partir de la décomposition de $\GL_2(\Zp)\Qp^\times \backslash \GL_2(\Qp) / \GL_2(\Zp)\Qp^\times$ et on peut montrer qu'elle est isomorphe à $k_E[T]$, où $T$ correspond à la double classe $\GL_2(\Zp)\Qp^\times \cdot \smat{p & 0 \\ 0 & 1} \cdot \GL_2(\Zp)$.

Si $\chi : \Qp^\times \to k_E^\times$ est un caractère lisse et si $\lambda \in k_E$, alors on pose
\[ \pi(r,\lambda,\chi) = \frac{ \ind_{\GL_2(\Zp)\Qp^\times}^{\GL_2(\Qp)} \Sym^r k_E^2}{T-\lambda} \otimes (\chi \circ \det). \]
C'est une représentation lisse de $\GL_2(\Qp)$, de caractère central $\omega^r \chi^2$.

\begin{theo}\label{bbl1}
Si $r \in \{0,\hdots,p-1\}$ et si $(r,\lambda) \notin \{(0,\pm 1), (p-1, \pm 1)\}$, alors la représentation $\pi(r,\lambda,\chi)$ est irréductible. 

Si $\lambda = \pm 1$, alors on a deux suites exactes:
\begin{gather*} 
0 \to \Sp \otimes (\chi\mu_\lambda \circ \det) \to \pi(0,\lambda,\chi) \to \chi\mu_\lambda \circ \det \to 0, \\
0 \to \chi\mu_\lambda \circ \det \to \pi(p-1,\lambda,\chi) \to \Sp \otimes (\chi\mu_\lambda \circ \det) \to 0,
\end{gather*}
où la représentation $\Sp$ (la \emph{spéciale}) ainsi définie est irréductible.
\end{theo}

Ce théorème est la réunion des résultats de \cite{BL1} et \cite{BL2} qui traitent le cas $\lambda \neq 0$ et des résultats de \cite{BR1} qui traite le cas $\lambda=0$ (les  représentations dites \emph{supersingulières}).

\begin{theo}\label{bbl2}
Les représentations $k_E$-linéaires lisses absolument irréductibles de $\GL_2(\Qp)$ admettant un caractère central sont les suivantes :
\begin{enumerate}
\item $\chi \circ \det$;
\item $\Sp \otimes (\chi \circ \det)$;
\item $\pi(r,\lambda,\chi)$ où $r \in \{0,\hdots,p-1\}$ et $(r,\lambda) \notin \{(0,\pm 1), (p-1, \pm 1)\}$.
\end{enumerate}
\end{theo}

Ce théorème est démontré dans \cite{BL1}, \cite{BL2} et \cite{BR1}. On constate alors que toutes les représentations lisses irréductibles de $\GL_2(\Qp)$ admettant un caractère central sont \emph{admissibles}, c'est-à-dire que si $\K$ est un sous-groupe ouvert compact de $\GL_2(\Qp)$, alors l'espace des vecteurs fixes par $\K$ est de dimension finie. 

Il existe une autre manière de construire des représentations lisses de $\GL_2(\Qp)$, en utilisant l'induction parabolique. On note $\B_2(\Qp)$ l'ensemble des matrices triangulaires supérieures de $\GL_2(\Qp)$. Si $\delta_1$ et $\delta_2$ sont deux caractères lisses de $\Qp^\times$, alors on note $\delta_1 \otimes \delta_2$ le caractère de $\B_2(\Qp)$ défini par $\smat{a & b \\ 0 & d} \mapsto \delta_1(a) \delta_2(d)$. La représentation 
\[ B(\delta_1,\delta_2) = \ind_{\B_2(\Qp)}^{\GL_2(\Qp)} (\delta_1 \otimes \omega^{-1} \delta_2) \]
est l'ensemble des fonctions $f : \GL_2(\Qp) \to k_E$ localement constantes et telles que $f(bg)= (\delta_1 \otimes \omega^{-1} \delta_2)(b) f(g)$ si $b \in \B_2(\Zp)$ et $g \in \GL_2(\Qp)$. Cet espace est muni de l'action de $\GL_2(\Qp)$ donnée par $(gf)(h)=f(hg)$. On n'obtient pas de nouvelles représentations, comme le précise le résultat suivant.

\begin{theo}\label{bl3}
Si $\lambda \in k_E^\times$ et $r \in \{0,\hdots,p-1\}$, alors les semi-simplifiées des représentations $B(\chi\mu_{1/\lambda} \otimes \chi\omega^{r+1} \mu_\lambda)$ et $\pi(r,\lambda,\chi)$ sont isomorphes.
\end{theo}

Ces isomorphismes sont démontrés dans \cite{BL1} et \cite{BL2}.

Notons que la représentation $B(1,\omega)$ se réalise (via l'identification entre $\Proj^1(\Qp)$ et $\GL_2(\Qp)/\B_2(\Qp)$) comme l'espace $\calC^0(\Proj^1(\Qp),k_E)$ des fonctions localement constantes sur $\Proj^1(\Qp)$, muni de l'action naturelle de $\GL_2(\Qp)$. On trouve donc par le théorème \ref{bbl1} que la spéciale s'identifie à $\Sp \simeq \calC^0(\Proj^1(\Qp),k_E) / \{\text{constantes}\}$.

\subsection{La correspondance semi-simple modulo $p$}
\label{corrmodp}

Toute représentation $k_E$-linéaire semi-simple de dimension $2$ de $\Gal(\Qpbar/\Qp)$ est (après extension éventuelle des scalaires), soit absolument irréductible et donc de la forme $\rho(r,\chi)$ par le théorème \ref{repskegal}, soit une somme de deux caractères et de la forme $(\omega^{r+1} \mu_\lambda \oplus \mu_{1/\lambda}) \otimes \chi$ avec $\lambda \in k_E^\times$ et $r \in \{0,\hdots,p-2\}$. 

La correspondance entre représentations $k_E$-linéaires semi-simples de dimension $2$ de $\Gal(\Qpbar/\Qp)$ et représentations $k_E$-linéaires semi-simples lisses de $\GL_2(\Qp)$, est définie par Breuil dans \cite{BR1} comme suit
\begin{align*}
\rho(r,\chi) & \leftrightarrow \pi(r,0,\chi), \\
(\omega^{r+1} \mu_\lambda \oplus \mu_{1/\lambda}) \otimes \chi &  \leftrightarrow \pi(r,\lambda,\chi)^{\ssf} \oplus \pi([p-3-r],1/\lambda,\omega^{r+1} \chi)^{\ssf},
\end{align*}
où $[p-3-r]$ est le représentant de $p-3-r \bmod{p-1}$ dans $\{0,\hdots,p-2\}$ et {\og $\ssf$\fg}~veut dire semi-simplifiée. Du côté $\GL_2(\Qp)$, les objets peuvent être de longueur $1$, $2$, $3$~ou~$4$.

\begin{theo}\label{cmbd}
La correspondance ci-dessus est bien définie.
\end{theo}

Il s'agit de vérifier que des paramètres $(r,\lambda,\chi)$ différents qui donnent des représentations isomorphes d'un côté donnent des représentations isomorphes de l'autre côté. On sait par exemple que les seuls entrelacements entre les $\rho(r,\chi)$ sont 
\[ \rho(r,\chi) = \rho(r,\chi\mu_{-1})= \rho(p-1-r,\chi\omega^r) =\rho(p-1-r,\chi\omega^r\mu_{-1}), \]
et que de même les seuls entrelacements entre les $\pi(r,0,\chi)$ sont 
\[ \pi(r,0,\chi) = \pi(r,0,\chi\mu_{-1})= \pi(p-1-r,0,\chi\omega^r) =\pi(p-1-r,0,\chi\omega^r\mu_{-1}). \]

La correspondance ci-dessus a été définie par Breuil, en se fondant sur les calculs de réduction modulo $p$ rappelés au \S \ref{prcons}. En étant plus soigneux, on peut aussi définir une correspondance sans semi-simplifier.

\section{Représentations en caractéristique $0$}
\label{repnul}

Dans cette section, nous rappelons la théorie de Fontaine pour les représentations $p$-adiques de $\Gal(\Qpbar/\Qp)$, et la théorie de Schneider et Teitelbaum pour les représentations de $\GL_2(\Qp)$. Ensuite, nous donnons les premiers exemples de la correspondance $p$-adique.

\subsection{Théorie de Hodge $p$-adique}
\label{thpad}

Il est facile de faire la liste des représentations $k_E$-linéaires de dimension $2$ de $\Gal(\Qpbar/\Qp)$, mais l'étude des représentations $E$-linéaires de ce groupe est plus compliquée car il y en a beaucoup. De manière plus précise, si l'on se donne une représentation $k_E$-linéaire $W$ de dimension $2$ de $\Gal(\Qpbar/\Qp)$, alors les résultats de \cite{MDG} montrent que l'ensemble des représentations $E$-linéaires $V$ qui admettent un $\OO_E$-réseau $T$ stable par $\Gal(\Qpbar/\Qp)$ avec $T/\MM_E T = W$ est en général l'ensemble des $\OO_E$-points d'un espace rigide de dimension $5$.

L'objet de la théorie de Hodge $p$-adique est de faire le tri dans ces représentations, et de décrire aussi explicitement que possible celles qui proviennent de la géométrie arithmétique. Pour cela, Fontaine a introduit dans \cite{FPP} un certain nombre d'anneaux $\bcris \subset \bst \subset \bdr$ qui sont des $\Qp$-algèbres topologiques munies d'une action de $\Gal(\Qpbar/\Qp)$. L'anneau $\bdr$ est muni d'une filtration et l'anneau $\bst$ est muni d'un frobenius $\phi$ et d'un opérateur de monodromie $N$ tels que $N \circ \phi = p \phi \circ N$; enfin, on a 
\mbox{$\bcris=\bst^{N=0}$.} 
Si $V$ est une représentation $E$-linéaire de $\Gal(\Qpbar/\Qp)$ et si 
\mbox{$* \in \{$cris, st, dR$\}$,} alors on pose $\dfont_*(V) = (\bfont_*  \otimes_{\Qp} V)^{G_K}$. On peut montrer que $\dfont_*(V)$ est un $E$-espace vectoriel de dimension $\leq \dim_E(V)$ et on dit que $V$ est \emph{cristalline} ou \emph{semi-stable} ou de \emph{de Rham} si l'on a égalité pour $*$ égal à cris ou st ou dR.

Le $E$-espace vectoriel $\ddr(V)$ est alors muni d'une filtration par des sous-espaces $E$-linéaires, l'espace $\dst(V) \subset \ddr(V)$ est un $(\phi,N)$-module filtré et $\dcris(V) = \dst(V)^{N=0}$. Si $D$ est un $(\phi,N)$-module filtré, alors on définit $t_N(D)$ comme étant la valuation $p$-adique de $\phi$ sur $\det(D)$ et $t_H(D)$ comme étant l'unique entier $h$ tel que $\Fil^h(\det(D))=\det(D)$ et $\Fil^{h+1}(\det(D))=\{0\}$. On dit que $D$ est admissible si $t_H(D)=t_N(D)$ et si $t_H(D') \leq t_N(D')$ pour tout sous-objet $D'$ de $D$.

\begin{theo}\label{cf}
Si $V$ est une représentation semi-stable de $\Gal(\Qpbar/\Qp)$, alors $\dst(V)$ est un $(\phi,N)$-module filtré admissible et le foncteur $V \mapsto \dst(V)$ donne une équivalence de catégories de la catégorie des représentations semi-stables vers la catégorie des $(\phi,N)$-modules filtrés admissibles.
\end{theo}

Le fait que $\dst(V)$ est admissible et que le foncteur $V \mapsto \dst(V)$ est pleinement fidèle est démontré dans \cite{FST}. Le fait que tout module admissible provient d'une représentation semi-stable est le résultat principal de \cite{CF}. L'intérêt de ce théorème est que pour se donner une représentation semi-stable, il suffit de se donner un $(\phi,N)$-module filtré admissible, qui est un objet tout à fait explicite. 

\subsection{Représentations admissibles de $\GL_2(\Qp)$}
\label{admgl}

Les représentations de $\GL_2(\Qp)$ qui nous intéressent sont les \emph{$\GL_2(\Qp)$-banach unitaires}, c'est-à-dire les espaces de Banach $\B$ munis d'une action continue de $\GL_2(\Qp)$ et dont la topologie est définie par une norme $\|{\cdot}\|$ telle que $\|g(v)\| = \|v\|$ quels que soient $g \in \GL_2(\Qp)$ et $v \in \B$. 

Si $\B$ est un tel objet, alors $\OO_{\B} = \{ b \in \B$ tels que $\|b\| \leq 1\}$ est stable par $\GL_2(\Qp)$ et donc $\overline{\B} = \OO_{\B} / \MM_E \OO_{\B}$ est une représentation $k_E$-linéaire lisse de $\GL_2(\Qp)$. On dit comme dans \cite{STI} qu'un $\GL_2(\Qp)$-banach unitaire $\B$ est \emph{admissible} si $\overline{\B}$ est admissible au sens usuel, c'est-à-dire que si $\K$ est un sous-groupe ouvert compact de $\GL_2(\Qp)$, alors $\overline{\B}^{\K}$ est de dimension finie. Si $\overline{\B}$ est de longueur finie, alors la semi-simplifiée de $\overline{\B}$ ne dépend pas du choix de la norme par le principe de Brauer-Nesbitt, et on l'appelle par abus de langage la réduction modulo $\MM_E$ de $\B$.

Si l'on prend par exemple $\B=\calC^0(\Proj^1(\Qp),E)$, alors $\B$ est un $\GL_2(\Qp)$-banach unitaire admissible. L'espace $\B$ contient de manière évidente les constantes et la \emph{Steinberg} est par définition $\St = \calC^0(\Proj^1(\Qp),E) / \{ \text{constantes} \}$. C'est une représentation irréductible admissible de $\GL_2(\Qp)$, dont la réduction modulo $\MM_E$ est la spéciale $\Sp$.

Si $\B$ est un $\GL_2(\Qp)$-banach unitaire et si $v \in \B$, alors on a une fonction 
\mbox{$\GL_2(\Qp) \to \B$} donnée par $g \mapsto g(v)$ et on dit que $v$ est \emph{localement analytique} si $g \mapsto g(v)$ l'est. On dit de même que $v$ est \emph{localement algébrique} si $g \mapsto g(v)$ est localement un polynôme en $a$, $b$, $c$, $d$ et $(ad-bc)^{-1}$ où $g = \smat{a & b \\ c & d}$, et que $v$ est \emph{localement constante} (c'est-à-dire lisse) si $g \mapsto g(v)$ l'est. On note $\B^{\an}$, $\B^{\alg}$ et $\B^{\lisse}$ les sous-espaces vectoriels correspondants de $\B$. En général, on peut très bien avoir $\B^{\alg} = \{ 0 \}$ et $\B^{\lisse} = \{0\}$, mais on a le résultat suivant qui est démontré dans \cite{ST1}.

\begin{theo}\label{st}
Si $\B$ est un $\GL_2(\Qp)$-banach unitaire admissible, alors $\B^{\an}$ est un sous-espace dense dans $\B$.
\end{theo}

Ce résultat est à comparer avec le théorème correspondant de Vignéras dans \cite{MBLP}, concernant les représentations $\ell$-adiques de groupes $p$-adiques avec $\ell \neq p$. Dans ce cas, les vecteurs lisses forment déjà un sous-espace dense.

\subsection{Premiers exemples}
\label{prcons}

Avant d'expliquer au \S\ref{constcorr} la construction générale de la correspondance de Langlands $p$-adique pour $\GL_2(\Qp)$, nous donnons ici les premiers exemples historiques, construits par Breuil dans \cite{BR2} et \cite{BRL}. 

Commençons par le cas cristallin. On se donne un entier $k \geq 2$ et un nombre $a_p \in \MM_E$. On associe à ces données le $\phi$-module filtré $D_{k,a_p} = E e_1 \oplus E e_2$ où
\[ \Mat(\phi) = \pmat{0 & -1 \\ p^{k-1} & a_p} 
\text{ et } \Fil^i D_{k,a_p} = 
\begin{cases} 
D_{k,a_p} & \text{ si $i \leq 0$,} \\
E e_1 & \text{ si $1 \leq i \leq k-1$,} \\
\{0\} & \text{ si $i \geq k$.}
\end{cases} \]
Ce $\phi$-module filtré est admissible, et il existe donc une représentation cristalline $V_{k,a_p}$ telle que $\dcris(V_{k,a_p}^*) = D_{k,a_p}$. Si $\chi : \Gal(\Qpbar/\Qp) \to \OO_E^\times$ est un caractère, alors on note $V_{k,a_p,\chi} = V_{k,a_p} \otimes \chi$.

\begin{theo}\label{vkap}
Toute représentation cristalline absolument irréductible de dimension $2$ de $\Gal(\Qpbar/\Qp)$ est de la forme $V_{k,a_p,\chi}$ avec $k \geq 2$, $a_p \in \MM_E$ et 
\mbox{$\chi : \Gal(\Qpbar/\Qp) \to \OO_E^\times$} un caractère cristallin. De plus, les seuls isomorphismes entre ces représentations sont donnés par $V_{k,a_p,\chi} = V_{k,-a_p,\chi\mu_{-1}}$.
\end{theo}

Si $k \geq 2$, soit $\Sym_0^{k-2} E^2$ la représentation $\Sym^{k-2} E^2$ de $\GL_2(\Zp)$, étendue à $\GL_2(\Zp)\Qp^\times$ en envoyant $p$ sur l'identité. On pose alors 
\[ \Pi_{k,a_p,\chi} = \frac{\ind_{\GL_2(\Zp)\Qp^\times}^{\GL_2(\Qp)} \Sym_0^{k-2} E^2}{T-a_p} \otimes (\chi \circ \det), \]
où $T$ est un opérateur de Hecke défini comme en \ref{gl2modp}. Cette représentation est alors localement algébrique, irréductible si $a_p \neq \pm(p^{k/2}+p^{k/2-1})$, et peut se réaliser comme le produit tensoriel d'une représentation algébrique par une représentation lisse
\[ \Pi_{k,a_p,\chi} = \Sym^{k-2} E^2 \otimes \ind_{\B_2(\Qp)}^{\GL_2(\Qp)} (\mu_{\lambda_1} \otimes \mu_{p \lambda_2^{-1}}) \otimes (\chi \circ \det), \]
où $\lambda_1$ et $\lambda_2$ sont les racines de $X^2-a_pX+p^{k-1}=0$. Enfin, les seuls isomorphismes entre ces représentations sont donnés par $\Pi_{k,a_p,\chi} = \Pi_{k,-a_p,\chi\mu_{-1}}$.

\begin{theo}\label{pikapcor}
La représentation $\Pi_{k,a_p,\chi}$ admet un réseau de type fini stable sous $\GL_2(\Qp)$, et son complété pour ce réseau est un $\GL_2(\Qp)$-banach unitaire admissible et topologiquement irréductible.
\end{theo}

Notons que deux réseaux de type fini sont commensurables et donnent donc des complétés isomorphes. Le complété de $\Pi_{k,a_p,\chi}$ est topologiquement irréductible y compris quand $\Pi_{k,a_p,\chi}$ est réductible (ce qui se produit pour $a_p = \pm(p^{k/2}+p^{k/2-1})$). Le théorème \ref{pikapcor} avait été conjecturé par Breuil, et démontré pour $k \leq 2p$ dans le \S 3.3 de \cite{BR2}. La démonstration consistait à trouver explicitement un réseau et à calculer sa réduction modulo $\MM_E$. Une démonstration générale est donnée dans \cite{LB14} (le cas $a_p = \pm 2 p^{(k-1)/2}$ pose un problème particulier, voir \cite{PCR}). 

\begin{theo}\label{pikapred}
La semi-simplifiée $\overline{V}_{k,a_p,\chi}$ de la réduction modulo $\MM_E$ de $V_{k,a_p,\chi}$ correspond (via la correspondance du \S\ref{corrmodp}) à la semi-simplifiée $\overline{\Pi}_{k,a_p,\chi}$ de la réduction modulo $\MM_E$ de $\Pi_{k,a_p,\chi}$.
\end{theo}

Ce résulat suivait des calculs de Breuil dans tous les cas où $\overline{V}_{k,a_p,\chi}$ était connu à l'époque (soit par la théorie de Fontaine-Laffaille de \cite{FL}, soit par des calculs informatiques explicites d'exemples comme dans \cite{SSB}), et renforçait l'idée que les $\Pi_{k,a_p,\chi}$ étaient les bons objets. Le théorème \ref{pikapred} a été démontré en toute généralité dans \cite{LB15} en s'appuyant sur les constructions de \cite{LB14} rappelées en \ref{triang}.

Breuil a par ailleurs donné dans \cite{BRL} une construction similaire pour des représen\-tations semi-stables non cristallines. Si $k > 2$ et $\Linv \in E$, soit $D_{k,\Linv}$ le $(\phi,N)$-module filtré $D_{k,\Linv} = E e_1 \oplus E e_2$ avec $\Mat(N) = \smat{0 & 0 \\ 1 & 0}$ et 
\[ \Mat(\phi) = \pmat{p^{k/2} & 0 \\ 0 & p^{k/2-1}} 
\text{ et } \Fil^i D_{k,\Linv} = 
\begin{cases} 
D_{k,\Linv} & \text{ si $i \leq 0$,} \\
E (e_1+ \Linv e_2) & \text{ si $1 \leq i \leq k-1$,} \\
\{0\} & \text{ si $i \geq k$.}
\end{cases} \]
Ce $(\phi,N)$-module filtré est admissible et il existe donc une représentation semi-stable $V_{k,\Linv}$ telle que $\dst(V_{k,\Linv}^*)=D_{k,\Linv}$.

Soit $\log_{\Linv}$ le logarithme $p$-adique normalisé par $\log_{\Linv}(p)=\Linv$ et $W(\Linv)$ la $E$-représen\-tation de dimension $2$ de $\B_2(\Qp)$ donnée par 
\[ \pmat{a & b \\ 0 & d} \mapsto \pmat{1 & \log_{\Linv}(a/d) \\ 0 & 1}, \]
ce qui fait que $W(\Linv)$ est une extension non-scindée de $E$ par $E$. Si $\delta_k : \B_2(\Qp) \to E^\times$ est le caractère qui à $\smat{a & b \\ 0 & d}$ associe $|ad|^{(k-2)/2} d^{k-2}$, alors on pose $W(k,\Linv)=W(\Linv) \otimes \delta_k$ et on a une suite exacte
\[ 0 \to \ind_{\B_2(\Qp)}^{\GL_2(\Qp)} \delta_k \to  \ind_{\B_2(\Qp)}^{\GL_2(\Qp)} W(k,\Linv) \xrightarrow{s} 
\ind_{\B_2(\Qp)}^{\GL_2(\Qp)} \delta_k \to 0. \]
Par ailleurs, on peut montrer que la représentation $\Sym^{k-2} E^2 \otimes |{\det}|^{(k-2)/2}$ est une sous-représentation de $\ind_{\B_2(\Qp)}^{\GL_2(\Qp)} \delta_k$ et on définit 
\[ \Sigma(k,\Linv) = \frac{s^{-1}(\Sym^{k-2} E^2 \otimes |{\det}|^{(k-2)/2})}{\Sym^{k-2} E^2 \otimes |{\det}|^{(k-2)/2}}. \]
Les analogues des théorèmes \ref{pikapcor} et \ref{pikapred} sont alors vrais : la représentation $\Sigma(k,\Linv)$ admet un réseau stable sous $\GL_2(\Qp)$, et le complété $B(k,\Linv)$ de $\Sigma(k,\Linv)$ est un $\GL_2(\Qp)$-banach unitaire admissible topologiquement irréductible. De plus, la semi-simplifiée de la réduction de $V_{k,\Linv}$ modulo $\MM_E$ correspond à la semi-simplifiée de celle de $B(k,\Linv)$ via la correspondance du \S\ref{corrmodp} (voir \cite{BMC} pour des cas particuliers, qui ont eux aussi renforcé l'idée que les $B(k,\Linv)$ étaient les bons objets; ce résultat est démontré en toute généralité dans \cite{LB15} en s'appuyant cette fois sur les constructions de \cite{CSP} rappelées en \ref{triang}).

\section{La série principale unitaire}
\label{constcorr}

Dans cette section, nous rappelons la théorie des $(\phi,\Gamma)$-modules de Fontaine, puis nous expliquons  son application à la construction de modèles de la restriction à $\B_2(\Qp)$ des représentations de $\GL_2(\Qp)$ associées aux représentations $p$-adiques {\og triangulines\fg}. 

\subsection{Les $(\phi,\Gamma)$-modules}
\label{pgmfont}

La théorie des $(\phi,\Gamma)$-modules de Fontaine, introduite dans \cite{F90}, permet de décrire toutes les représentations $p$-adiques de $\Gal(\Qpbar/\Qp)$ au moyen de modules sur des anneaux de séries munis de certains opérateurs. 

Soit $\calE^\dag$ le corps des séries $f(X)=\sum_{n \in \ZZ} a_n X^n$ où $a_n \in E$, la suite $\{a_n\}_{n \in \ZZ}$ est bornée et il existe $\rho(f)<1$ tel que $f(X)$ converge sur $\rho(f) \leq |X| < 1$. Cet anneau peut être muni de deux topologies : d'une part la topologie $p$-adique (la norme de Gauss), et d'autre part la topologie LF consistant à mettre sur chaque $\calE^{\dag,\rho}$ (les $f(X)$ qui convergent sur $\rho \leq |X| < 1$) la topologie de Fréchet de la convergence uniforme sur les couronnes du type $\rho \leq |X| \leq \sigma$ pour $\sigma<1$. 

Le complété de $\calE^\dag$ pour la topologie $p$-adique est le corps local $\calE$ des séries 
$f(X)=\sum_{n \in \ZZ} a_n X^n$ où $a_n \in E$, la suite $\{a_n\}_{n \in \ZZ}$ est bornée et $a_{-n} \to 0$ quand $n \to +\infty$. On note $\OO_{\calE}^\dag$ et $\OO_{\calE}$ les anneaux des entiers de $\calE^\dag$ et de $\calE$ pour la norme de Gauss. 

Le complété de $\calE^\dag$ pour la topologie LF est \emph{l'anneau de Robba} $\calR$ des séries 
$f(X)=\sum_{n \in \ZZ} a_n X^n$ où $a_n \in E$ et il existe $\rho(f)<1$ tel que $f(X)$ converge sur $\rho(f) \leq |X| < 1$.

Soit $\Gamma$ un groupe isomorphe à $\Zp^\times$, dont on note $[a]$ l'élément correspondant à 
\mbox{$a \in \Zp^\times$.} Tous les anneaux ci-dessus sont munis d'un frobenius $\phi$ défini 
par $\phi(f)(X)=$\break $f((1+X)^p-1)$ et d'une action de $\Gamma$ donnée par $([a]f)(X)=f((1+X)^a-1)$. 

\begin{defi}
Si $A$ est l'un des anneaux $\calE^\dag$ ou $\calE$ ou $\calR$, alors un \emph{$(\phi,\Gamma)$-module sur $A$} est un $A$-module libre de rang fini $d$, muni d'un frobenius semi-linéaire $\phi$ tel que $\Mat(\phi) \in \GL_d(A)$ et d'une action semi-linéaire et continue de $\Gamma$ qui commute à $\phi$. 

On dit qu'un $(\phi,\Gamma)$-module sur $A$ est \emph{étale} s'il en existe une base dans laquelle $\Mat(\phi) \in \GL_d(\OO_{\calE}^\dag)$ (si $A$ est $\calE^\dag$ ou $\calR$) ou dans laquelle $\Mat(\phi) \in \GL_d(\OO_{\calE})$ (si $A = \calE$).
\end{defi}

Si $\dfont^\dag$ est un $(\phi,\Gamma)$-module sur $\calE^\dag$, alors $\dfont = \calE \otimes_{\calE^\dag} \dfont^\dag$ est un $(\phi,\Gamma)$-module sur $\calE$ et $\drig = \calR \otimes_{\calE^\dag} \dfont^\dag$ est un $(\phi,\Gamma)$-module sur $\calR$. De plus, si $\dfont^\dag$ est étale, alors $\dfont$ et $\drig$ le sont aussi.

\begin{theo}\label{ccked}
Les foncteurs $\dfont^\dag \mapsto \dfont$ et $\dfont^\dag \mapsto \drig$ sont des équivalences de catégories, de la catégorie des $(\phi,\Gamma)$-modules étales sur $\calE^\dag$, vers la catégorie des \mbox{$(\phi,\Gamma)$-modules} étales sur $\calE$ et sur $\calR$.
\end{theo}

Le fait que $\dfont^\dag \mapsto \dfont$ est une équivalence de catégories (la {\og surconvergence\fg} des \mbox{$(\phi,\Gamma)$-modules} sur $\calE$) est le résultat principal de \cite{CC98}. Le fait que $\dfont^\dag \mapsto \drig$ est une équivalence de catégories est démontré dans \cite{KLMT}.

Il existe un anneau ${\calE}^{\nr}$ défini par exemple dans \cite{F90} (c'est le complété de 
\mbox{l'extension} maximale non-ramifiée de $\calE$), qui est muni d'un frobenius $\phi$ et d'une action de $\Gal(\Qpbar/\Qp)$ et qui contient $\calE$, ce qui fait que si $\dfont$ est un $(\phi,\Gamma)$-module sur $\calE$, alors $V(\dfont) = ({\calE}^{\nr} \otimes_{\calE} \dfont)^{\phi=1}$ est un $E$-espace vectoriel, muni de l'action de $\Gal(\Qpbar/\Qp)$ donnée par $g(x \otimes d) = g(x) \otimes [\chi_{\mathrm{cycl}}(g)](d)$.

\begin{theo}\label{dagrep}
Si $\dfont$ est un $(\phi,\Gamma)$-module étale de dimension $d$ sur $\calE$, alors $V(\dfont)$ est une représentation $E$-linéaire de dimension $d$ de $\Gal(\Qpbar/\Qp)$ et le foncteur qui en résulte, de la catégorie des $(\phi,\Gamma)$-modules étales sur $\calE$ vers la catégorie des représentations $E$-linéaires de $\Gal(\Qpbar/\Qp)$, est une équivalence de catégories.
\end{theo}

Afin de faire le lien entre les $(\phi,\Gamma)$-modules et les représentations de $\GL_2(\Qp)$, il faut construire un certain opérateur noté $\psi$. Si $A$  est l'un des anneaux $\calE^\dag$ ou $\calE$ ou $\calR$, alors $A$ est un $\phi(A)$-module libre de rang $p$ engendré par $\{(1+X)^i\}_{0 \leq i \leq p-1}$. Si $f \in A$, on peut donc écrire $f = \sum_{i=0}^{p-1} \phi(f_i) (1+X)^i$ et on pose $\psi(f)=f_0$. Si $\dfont$ est un $(\phi,\Gamma)$-module sur $A$, alors il en existe une base de la forme $\{\phi(e_i)\}_{1 \leq i \leq d}$ et si $y \in \dfont$, on peut donc écrire $y = \sum_{i=1}^{d} y_i \phi(e_i)$ et on pose $\psi(y)=\sum_{i=1}^d \psi(y_i)e_i$.

\begin{prop}\label{psid}
L'opérateur $\psi$ ainsi défini ne dépend pas des choix, commute à l'action de $\Gamma$, et vérifie $\psi(\phi(f)y))=f\psi(y)$ et $\psi(f \phi(y))=\psi(f)y$ si $f \in A$.
\end{prop}

Le point de départ de la construction de la correspondance de Langlands $p$-adique pour $\GL_2(\Qp)$ en utilisant les $(\phi,\Gamma)$-modules est le suivant. Si $\dfont$ est un $(\phi,\Gamma)$-module sur $\calE$ et si $\delta:\Qp^\times \to E^\times$ est un caractère continu, notons $\dfont \boxtimes_\delta \Qp$ l'ensemble des suites $\{x^{(n)}\}_{n \in \ZZ}$ d'éléments de $\dfont$ telles que $\psi(x^{(n+1)})=x^{(n)}$ pour tout $n$. On munit $\dfont \boxtimes_\delta \Qp$ d'une action de $\B_2(\Qp)$ en décidant que si $x \in \dfont \boxtimes_\delta \Qp$, alors 
\begin{enumerate}
\item $g(x)^{(n)} = \delta(a) \cdot x^{(n)}$ si $g=\smat{a & 0 \\ 0 & a}$ avec $a \in \Qp^\times$;
\item $g(x)^{(n)} = [a](x^{(n)})$ si $g=\smat{a & 0 \\ 0 & 1}$ avec $a \in \Zp^\times$;
\item $g(x)^{(n)} = x^{(n+k)}$ si $g=\smat{p^k & 0 \\ 0 & 1}$ avec $k \in \ZZ$;
\item $g(x)^{(n)} = (1+X)^{c p^n} \cdot x^{(n)}$ si $g=\smat{1 & c \\ 0 & 1}$ et $c p^n \in \Zp$.
\end{enumerate}

Cette définition est donnée dans \cite{CPG}. L'idée est alors d'étendre cette action de $\B_2(\Qp)$ à une action de $\GL_2(\Qp)$. Cette stratégie marche particulièrement bien pour les représentations {\og triangulines \fg} de dimension $2$, que nous étudions dans le \S\ref{triang}.

\subsection{Pentes des $\phi$-modules sur l'anneau de Robba}
\label{kedslopes}

Un ingrédient important de la construction de la {\og série principale unitaire\fg} du \S \ref{triang} est la théorie des pentes de frobenius pour les $\phi$-modules sur $\calR$, théorie due à Kedlaya et développée dans \cite{KLMT}. On dit qu'un $\phi$-module sur $\calR$ est \emph{pur de pente $a/h$} s'il en existe une base dans laquelle $\Mat(p^{-a} \phi^h) \in \GL_d(\OO_{\calE}^\dag)$ (par exemple, être étale est équivalent à être pur de pente nulle). Un module pur d'une certaine pente est dit \emph{isocline}. Le résultat principal de la théorie des pentes est le théorème 6.10 de \cite{KLMT}.

\begin{theo}\label{kedsl}
Si $\dfont$ est un $\phi$-module sur $\calR$, alors il admet une unique filtration $\{0\} = \dfont_0 \subset \dfont_1 \subset \cdots \subset \dfont_\ell = \dfont$ par des sous-$\phi$-modules saturés, vérifiant
\begin{enumerate}
\item pour tout $i \geq 1$, le module $\dfont_i / \dfont_{i-1}$ est isocline;
\item si $s_i$ est la pente de $\dfont_i / \dfont_{i-1}$, alors $s_1 < s_2 < \cdots < s_\ell$.
\end{enumerate}
\end{theo}

Comme la filtration est unique, si $\dfont$ est en plus un $(\phi,\Gamma)$-module, alors les $\dfont_i$ sont eux aussi des $(\phi,\Gamma)$-modules.

Un point délicat mais crucial de la théorie des pentes est qu'un $\phi$-module sur $\calR$ qui est pur de pente $s$ n'admet pas de sous-objet de pente $<s$ par le théorème \ref{kedsl}, mais peut très bien admettre des sous-objets saturés de pente $> s$.

\subsection{La série principale unitaire}
\label{triang}

Si $V$ est une représentation $p$-adique, alors en combinant les théorèmes \ref{dagrep} et \ref{ccked}, on peut lui associer le $(\phi,\Gamma)$-module étale $\drig(V)$ sur $\calR$, qui est un objet de la catégorie de tous les $(\phi,\Gamma)$-modules sur $\calR$.

\begin{defi}\label{deftri}
On dit qu'une représentation $p$-adique $V$ est \emph{trianguline} si $\drig(V)$ est une extension successive de $(\phi,\Gamma)$-modules de rang $1$ sur $\calR$.
\end{defi}

En utilisant les résultats de \cite{LB2} qui font le lien entre la théorie de Hodge $p$-adique et la théorie des $(\phi,\Gamma)$-modules, on peut montrer le résultat suivant.

\begin{theo}
Les représentations semi-stables sont triangulines.
\end{theo}

Si $\delta : \Qp^\times \to E^\times$ est un caractère continu, alors $w(\delta)=\log_p \delta(u)/\log_p u$ ne dépend pas de $u \in 1+p\Zp$ et est appelé le \emph{poids de $\delta$}. La \emph{pente de $\delta$} est $u(\delta)=\vp(\delta(p))$. 

On définit $\calR(\delta)$ comme étant le $(\phi,\Gamma)$-module de rang $1$ engendré par $e_\delta$ avec 
\mbox{$\phi(e_\delta)=\delta(p)e_\delta$} et $[a](e_\delta) = \delta(a)e_\delta$. La pente de $\calR(\delta)$ au sens du \S \ref{kedslopes} est alors bien~$u(\delta)$.

\begin{theo}\label{pgr1}
Tout $(\phi,\Gamma)$-module de rang $1$ sur $\calR$ est isomorphe à $\calR(\delta)$ pour un caractère $\delta : \Qp^\times \to E^\times$.
\end{theo}

Si $V$ est une représentation trianguline de dimension $2$, alors $\drig(V)$ est une extension de deux $(\phi,\Gamma)$-modules de rang $1$ et on a donc une suite exacte
\[ 0 \to \calR(\delta_1) \to \drig(V) \to \calR(\delta_2) \to 0. \]
Le fait que $\drig(V)$ est étale force les relations $u(\delta_1) + u(\delta_2)=0$ et (à cause du théorème \ref{kedsl}) $u(\delta_1) \geq 0$. Si $u(\delta_1)=u(\delta_2)=0$, alors $\calR(\delta_1)$ et $\calR(\delta_2)$ sont étales et $V$ elle-même est extension de deux représentations. 

\begin{theo}\label{extpgr}
Si $\delta_1$ et $\delta_2 : \Qp^\times \to E^\times$ sont deux caractères continus, alors $\Ext^1(\calR(\delta_2), \calR(\delta_1))$ est un $E$-espace vectoriel de dimension $1$, sauf si $\delta_1 \delta_2^{-1}$ est de la forme $x^{-i}$ avec $i \geq 0$, ou de la forme $|x|x^i$ avec $i \geq 1$; dans ces deux cas, $\Ext^1(\calR(\delta_2), \calR(\delta_1))$ est de dimension $2$.
\end{theo}

Ce théorème est démontré dans \cite{CTR}. On note alors $\calS$ l'espace $\calS=\{ (\delta_1,\delta_2,\Linv) \}$ où $\Linv = \infty$ si   $\delta_1 \delta_2^{-1}$ n'est pas de la forme $x^{-i}$ avec $i \geq 0$, ni de la forme $|x|x^i$ avec $i \geq 1$, et $\Linv \in \Proj^1(E)$ sinon. Une version constructive du théorème \ref{extpgr} ci-dessus permet d'associer à tout $s \in \calS$ une extension non-triviale $\drig(s)$ de $\calR(\delta_2)$ par $\calR(\delta_1)$, et réciproquement.

Si $s \in \calS$, alors on pose $w(s)=w(\delta_1)-w(\delta_2)$. On définit $\calS_*$ comme l'ensemble des $s \in \calS$ tels que $u(\delta_1)+u(\delta_2)=0$ et $u(\delta_1)>0$ et on pose alors $u(s)=u(\delta_1)$ si $s \in \calS_*$. On définit les ensembles {\og cristallins\fg}, {\og semi-stables\fg} et {\og non-géométriques\fg} de paramètres.
\begin{enumerate}
\item $\calS^{\cris}_* = \{ s \in \calS_*$ tels que $w(s) \geq 1$ et $u(s)<w(s)$ et $\Linv=\infty \}$;
\item $\calS^{\st}_* = \{ s \in \calS_*$ tels que $w(s) \geq 1$ et $u(s)<w(s)$ et $\Linv \neq \infty \}$;
\item $\calS^{\ngeo}_* = \{ s \in \calS_*$ tels que $w(s)$ n'est pas un entier $\geq 1 \}$;
\item $\calS_{\irr} = \calS^{\cris}_* \sqcup \calS^{\st}_* \sqcup \calS^{\ngeo}_*$.
\end{enumerate}

\begin{theo}\label{dsirr}
Si $s \in \calS_{\irr}$, alors $\drig(s)$ est étale et irréductible. Si $V(s)$ est la représentation trianguline associée, alors $V(s)=V(s')$ si et seulement si $s \in \calS^{\cris}_*$ et $s'=(x^{w(s)} \delta_2, x^{-w(s)} \delta_1,\infty)$. 
\end{theo}

Toutes les représentations triangulines absolument irréductibles s'obtiennent ainsi (quitte à étendre les scalaires) et $V(s)$ devient cristalline (ou semi-stable) sur une extension abélienne de $\Qp$ après torsion éventuelle par un caractère si $s \in \calS^{\cris}_*$ (ou si $s \in \calS^{\st}_*$), tandis qu'elle n'est pas de de Rham si $s \in \calS^{\ngeo}_*$. Un ingrédient important de la démonstration du théorème \ref{dsirr} donnée dans \cite{CTR} est la théorie des pentes des $\phi$-modules sur $\calR$ rappelée au \S \ref{kedslopes}. 

Nous expliquons à présent la construction des représentations $\Pi(s)$ de $\GL_2(\Qp)$ associées à $s \in \calS_{\irr}$. On note $\log_{\Linv}$ le logarithme normalisé par $\log_{\Linv}(p)=\Linv$ (si $\Linv=\infty$, on pose $\log_\infty=\vp$) et si $s \in \calS$, on note $\delta_s$ le caractère $(x|x|)^{-1} \delta_1 \delta_2^{-1}$. Si $s \in \calS_{\irr}$ alors on ne peut avoir $\Linv \neq \infty$ que si $\delta_s$ est de la forme $x^i$ avec $i \geq 0$. On peut définir la notion de fonction de classe $\calC^u$ pour $u \in \RR_{\geq 0}$ (voir \cite{CFV}), généralisant le cas $u \in \ZZ_{\geq 0}$. On note $\B(s)$ l'ensemble des fonctions $f : \Qp \to E$ qui sont de classe $\calC^{u(s)}$ et telles que $x \mapsto \delta_s(x) f(1/x)$ se prolonge en $0$ en une fonction de classe $\calC^{u(s)}$. L'espace $\B(s)$ est alors muni d'une action de $\GL_2(\Qp)$ donnée par la formule suivante
\[ \left[ \pmat{a & b \\ c & d} \cdot f \right] (y) = (x|x|\delta_1^{-1})(ad-bc) \cdot \delta_s(cy+d) \cdot f\left( \frac{ay+b}{cy+d}\right). \]

L'espace $M(s)$ est défini par
\begin{enumerate}
\item si $\delta_s$ n'est pas de la forme $x^i$ avec $i \geq 0$, alors $M(s)$ est l'espace engendré par $1$ et par les $y \mapsto \delta_s(y-a)$ avec $a \in \Qp$;
\item si $\delta_s$ est de la forme $x^i$ avec $i \geq 0$, alors $M(s)$ est l'intersection de $\B(s)$ et de l'espace engendré par $y \mapsto \delta_s(y-a)$ et $y \mapsto \delta_s(y-a)\log_{\Linv}(y-a)$ avec $a \in \Qp$.
\end{enumerate}
On pose enfin $\Pi(s)=\B(s)/\widehat{M}(s)$ où $\widehat{M}(s)$ est l'adhérence de $M(s)$ dans $\B(s)$.

\begin{theo}\label{pispikap}
Si $V(s)=V_{k,a_p,\chi}$, alors on a un morphisme $\GL_2(\Qp)$-équivariant $\Pi_{k,a_p,\chi} \to \Pi(s)$ dont l'image est dense, et si $V(s)=V_{k,\Linv}$, alors on a un morphisme $\GL_2(\Qp)$-équivariant $\Sigma(k,\Linv) \to \Pi(s)$ dont l'image est dense.
\end{theo}

Ce théorème n'exclut pas a priori que $\Pi(s)$ soit nul. On a cependant le résultat suivant, qui implique alors le théorème \ref{pikapcor} et son analogue semi-stable. On note $\delta$ le caractère $(x|x|)^{-1} \delta_1 \delta_2$.

\begin{theo}\label{spu}
Si $s \in \calS_{\irr}$, alors on a un isomorphisme de représentations de $\B_2(\Qp)$ entre $\Pi(s)^* \otimes \delta$ et l'ensemble des suites bornées de $\dfont(V(s)) \boxtimes_\delta \Qp$.
\end{theo}

Ce théorème avait tout d'abord été montré par Colmez pour $s \in\calS^{\st}_*$ puis par Breuil et moi-même pour $s \in \calS^{\cris}_*$ et enfin par Colmez pour $s \in \calS^{\ngeo}_*$ (voir \cite{LB14} et \cite{CSP}). Si $s \in \calS^{\cris}_*$ et $s'=(x^{w(s)} \delta_2, x^{-w(s)} \delta_1,\infty)$, alors $\Pi(s) = \Pi(s')$ et on a donc deux manières de construire cet espace. L'entrelacement entre $\Pi(s)$ et $\Pi(s')$ peut alors s'interpréter en termes de la filtration sur $\dcris(V(s))$.

\begin{coro}
Le $\GL_2(\Qp)$-banach unitaire $\Pi(s)$ est non-nul, topologiquement irréductible et admissible.
\end{coro}

Le théorème \ref{spu} implique que $\Pi(s)$ est non-nul car si $y \in \dfont(V(s))^{\psi=1}$, alors la suite constante de terme $y$ appartient à $\dfont(V(s)) \boxtimes_\delta \Qp$ et on peut montrer que $\dfont^{\psi=1} \neq 0$ pour tout $(\phi,\Gamma)$-module étale $\dfont$. Les deux autres propriétés se déduisent de même de la théorie des $(\phi,\Gamma)$-modules.

\begin{rema}\label{dsh}
On peut montrer que si $\dfont$ est un $(\phi,\Gamma)$-module étale, alors il existe un $\OO_E\dcroc{X}$-module {\og petit\fg} $\dsharp \subset \dfont$ tel que, si $x \in \dfont \boxtimes_\delta \Qp$ est une suite bornée, alors $x^{(n)} \in \dsharp$ pour tout $n$. Si $\dfont=\calE$, alors on peut prendre $\dsharp=X^{-1}\OO_E\dcroc{X} [1/p]$.
\end{rema}

\section{La correspondance pour $\GL_2(\Qp)$}
\label{constgl}

Cette section contient la construction générale de la correspondance, ainsi que plusieurs de ses propriétés.

\subsection{Les foncteurs de Colmez}
\label{foncolm}

Inspiré par ses constructions rappelées au \S\ref{triang}, Colmez a construit dans \cite{CGL} deux foncteurs $\dbf(\cdot)$ et $\Pibf(\cdot)$. Le foncteur $\dbf(\cdot)$ associe à un 
\mbox{$\GL_2(\Qp)$-banach} unitaire admissible $\Pi$ un $(\phi,\Gamma)$-module étale $\dbf(\Pi)$ sur $\calE$. Le foncteur $\Pibf(\cdot)$ associe à un \mbox{$(\phi,\Gamma)$-module} étale $\dfont$ absolument irréductible et de dimension $2$ sur $\calE$ un 
\mbox{$\GL_2(\Qp)$-banach} unitaire absolument irréductible et admissible $\Pibf(\dfont)$.

\textit{Construction de $\Pi \mapsto \dbf(\Pi)$} : si $\Pi$ est un $\GL_2(\Qp)$-banach unitaire admissible, alors on définit un certain sous-espace $W$ de $\Pi$ stable sous l'action du monoïde $P = \smat{\Zp\setminus\{0\} & \Zp \\ 0 & 1}$ ce qui fait que si $D_W$ est le dual de $W$, alors $D_W$ est une représentation de $P$. On peut donc munir $D_W$ d'une structure de $\OO_E\dcroc{X}$-module par $(1+X)^z \cdot v = \smat{1 & z \\ 0 & 1}v$, d'une action de $\Gamma$ par $[a](v) = \smat{a & 0 \\ 0 & 1}v$ et d'un frobenius $\phi$ par $\phi(v) = \smat{p & 0 \\ 0 & 1}v$. On pose alors $\dbf(\Pi) = \calE \otimes_{\OO_E\dcroc{X}} D_W$ et on vérifie que $\dbf(\Pi)$ est un $(\phi,\Gamma)$-module étale sur $\calE$.

\textit{Construction de $\dfont \mapsto \Pibf(\dfont)$} : si $\dfont$ est un $(\phi,\Gamma)$-module étale sur $\calE$ et si $\delta$ est un caractère, alors $\dfont \boxtimes_\delta \Qp$ est l'espace dont on a rappelé la définition à la fin de \ref{pgmfont}. Colmez a généralisé cette construction en définissant un faisceau $U \mapsto \dfont \boxtimes_\delta U$, avec des applications $\Res_U : \dfont \boxtimes_\delta V \to \dfont \boxtimes_\delta U$ si $U \subset V$ sont des ouverts de $\Qp$. On a par exemple $\dfont \boxtimes_\delta \Zp = \dfont$ et $\Res_{\Zp}$ est l'application $x = \{x^{(n)} \}_{n \in \ZZ} \mapsto x^{(0)}$. On a aussi $\dfont \boxtimes_\delta \Zp^\times = \dfont^{\psi=0}$, l'application $\Res_{\Zp^\times} : \dfont \boxtimes_\delta \Zp \to \dfont \boxtimes_\delta \Zp^\times$ étant donnée par $1-\phi\psi$. Colmez a alors défini par une formule explicite assez compliquée une application $w_\delta : \dfont \boxtimes_\delta \Zp^\times \to \dfont \boxtimes_\delta \Zp^\times$ qui correspond moralement à l'application $x \mapsto 1/x$ de $\Zp^\times \to \Zp^\times$. Cette application lui permet de définir $\dfont \boxtimes_\delta \Proj^1$ comme l'ensemble des $(x_1,x_2)$ avec $x_1,x_2 \in \dfont \boxtimes_\delta \Zp$ vérifiant $\Res_{\Zp^\times}(x_1) = w_\delta (\Res_{\Zp^\times} (x_2))$. On dispose alors d'une application $\Res_{\Qp} : \dfont \boxtimes_\delta \Proj^1 \to \dfont \boxtimes_\delta \Qp$ et l'action de $\B_2(\Qp)$ sur $\dfont \boxtimes_\delta \Qp$ s'étend assez naturellement en une action de $\GL_2(\Qp)$ sur $\dfont \boxtimes_\delta \Proj^1$ (on a par exemple $\smat{0 &  1 \\ 1 & 0}(x_1,x_2)=(x_2,x_1)$). 

Ces constructions marchent quelle que soit la dimension de $\dfont$, mais $\dfont \boxtimes_\delta \Proj^1$ n'est pas le genre d'objet que l'on cherche à construire (il est trop gros). Notons $(\dfont \boxtimes_\delta \Proj^1)^{\bor}$ le sous-module des $x \in \dfont \boxtimes_\delta \Proj^1$ tels que $\Res_{\Qp}(x)$ est une suite bornée de $\dfont \boxtimes_\delta \Qp$.

\begin{theo}\label{ccdf}
Si $\dfont$ est un $(\phi,\Gamma)$-module étale absolument irréductible et de dimension $2$ sur $\calE$ et si $\delta$ est le caractère $(x|x|)^{-1}\det(\dfont)$ alors
\begin{enumerate}
\item $(\dfont \boxtimes_\delta \Proj^1)^{\bor}$ est stable sous l'action de $\GL_2(\Qp)$;
\item si $\Pibf(\dfont) = \dfont \boxtimes_\delta \Proj^1 / (\dfont \boxtimes_\delta \Proj^1)^{\bor}$, alors $(\dfont \boxtimes_\delta \Proj^1)^{\bor}$ est isomorphe à $\Pibf(\dfont)^* \otimes \delta$;
\item on a un isomorphisme $\dbf(\Pibf(\dfont)) = \dfont \otimes \delta^{-1}$.
\end{enumerate}
\end{theo}

Si $\dfont$ n'est pas de dimension $2$, alors $(\dfont \boxtimes_\delta \Proj^1)^{\bor}$ n'est en général pas stable par $\GL_2(\Qp)$ et le théorème est donc spécifique à la dimension $2$. De plus, la démonstration du (1) est assez détournée, puisqu'elle est fondée sur le fait que le (1) est vrai pour les représentations triangulines par le théorème \ref{spu}, et que comme les représentations triangulines forment un sous-ensemble Zariski-dense de toutes les représentations \mbox{$p$-adiques,} le résultat s'étend par continuité (méthode suggérée par Kisin).

En utilisant l'équivalence de catégories entre représentations $p$-adiques et 
$(\phi,\Gamma)$-modules étales, on peut donc associer à toute représentation $E$-linéaire absolument irréductible de dimension $2$ de $\Gal(\Qpbar/\Qp)$ un $\GL_2(\Qp)$-banach unitaire absolument irréductible et admissible $\Pibf(V)$. Le résultat suivant de \cite{PMF} (valable si $p \geq 5$) nous dit quels $\GL_2(\Qp)$-banach on obtient de cette manière.

\begin{theo}\label{paskim}
Un $\GL_2(\Qp)$-banach unitaire absolument irréductible et admissible $\Pi$ est de la forme $\Pibf(V)$ avec $V$ absolument irréductible de dimension $2$ si et seulement si $\Pi$ n'est pas un sous-quotient d'une induite d'un caractère unitaire de $\B_2(\Qp)$.
\end{theo}

\begin{coro}\label{corrbij}
La correspondance de Langlands $p$-adique pour $\GL_2(\Qp)$ donne une bijection entre les deux ensembles suivants de $E$-représentations
\begin{enumerate}
\item les représentations absolument irréductibles de dimension $2$ de $\Gal(\Qpbar/\Qp)$
\item les $\GL_2(\Qp)$-banach unitaires absolument irréductibles et admissibles qui ne sont pas un sous-quotient d'une induite d'un caractère unitaire de $\B_2(\Qp)$.
\end{enumerate}
\end{coro}

Le résultat principal de \cite{KAC} nous dit par ailleurs que la plupart des représentations $E$-linéaires réductibles de dimension $2$ de $\Gal(\Qpbar/\Qp)$ sont dans l'image du foncteur $\dbf(\cdot)$ composé avec l'équivalence de Fontaine (l'idée est là aussi de se ramener au cas des triangulines par un argument de continuité). On peut alors étendre la correspondance aux représentations réductibles de dimension $2$ de $\Gal(\Qpbar/\Qp)$.

\subsection{Propriétés de la correspondance}
\label{props}

La correspondance de Langlands $p$-adique pour $\GL_2(\Qp)$ jouit d'un certain nombre de propriétés, qui ont d'ailleurs guidé sa construction.

Tout d'abord, elle est compatible à la correspondance en caractéristique $p$ du \S\ref{corrmodp}, par réduction modulo $\MM_E$. En effet, dans \cite{LB15} il est démontré que si $\Pi(W)$ est la représentation de $\GL_2(\Qp)$ associée à une représentation $k_E$-linéaire $W$ de $\Gal(\Qpbar/\Qp)$ par la correspondance du \S\ref{corrmodp}, alors la restriction à $\B_2(\Qp)$ de $\Pi(W)^* \otimes \delta$ est isomorphe à l'ensemble des suites bornées de $\dfont(W) \boxtimes_\delta \Qp$ après semi-simplification. En réduisant modulo $\MM_E$ l'isomorphisme entre $\Pibf(\dfont)^* \otimes \delta$ et $(\dfont \boxtimes_\delta \Proj^1)^{\bor}$ du théorème \ref{ccdf}, on obtient le résultat suivant.

\begin{theo}\label{compred}
Si $V$ est une représentation $E$-linéaire absolument irréductible de dimension $2$ de $\Gal(\Qpbar/\Qp)$, alors $\overline{V}$ correspond à $\overline{\Pibf}(V)$ par la correspondance du \S \ref{corrmodp}.
\end{theo}

Ensuite, la correspondance permet de retrouver la correspondance de Langlands locale {\og classique\fg} démontrée par Harris et Taylor dans \cite{HTLL} et par Henniart dans \cite{HSC} (mais construite bien avant par Tunnell dans \cite{TLL} pour $\GL_2(\Qp)$). Cette correspondance est une bijection entre des représentations du groupe de Weil-Deligne de $\Qp$ et des représentations lisses admissibles de $\GL_2(\Qp)$. Si $D$ est une représentation de Weil-Deligne de $\Qp$, on note $\LL(D)$ la représentation de $\GL_2(\Qp)$ associée. Si $D$ est de plus munie d'une filtration de poids $a<b$, alors on note $\AL(D)$ la représentation algébrique $\Sym^{b-a-1} E^2 \otimes \det^a$. Si $V$ est une représentation de $\Gal(\Qpbar/\Qp)$ qui est potentiellement semi-stable, alors une généralisation des constructions du \S\ref{thpad} permet de lui associer un $(\phi,N,G_{\Qp})$-module filtré $\dpst(V)$ et donc comme dans \cite{FWD} une représentation de Weil-Deligne et un module filtré. On note $\LL(V)$ et $\AL(V)$ les représentations de $\GL_2(\Qp)$ que l'on en déduit. Rappelons que si $\Pi$ est un $\GL_2(\Qp)$-banach unitaire, alors $\Pi^{\alg}$ a été défini au \S\ref{admgl}.

\begin{theo}\label{localg}
Si $V$ est une représentation $E$-linéaire absolument irréductible de dimension $2$ de $\Gal(\Qpbar/\Qp)$, alors $\Pibf(V)^{\alg} \neq \{0\}$ si et seulement si $V$ est potentiellement semi-stable à poids distincts $a<b$. Dans ce cas, on a $\Pibf(V)^{\alg} = \AL(V) \otimes \LL(V)$.
\end{theo}

Ce théorème est montré dans \cite{CGL} (en utilisant les résultats de \cite{EFC} si $V$ n'est pas trianguline) et nous ramène aux premiers exemples de la correspondance puisque les espaces $\Pi_{k,a_p,\chi}$ du \S \ref{prcons} pouvaient aussi être définis (si $a_p \neq \pm(p^{k/2}+p^{k/2-1})$) par $\Pi_{k,a_p,\chi} = \AL(V_{k,a_p,\chi}) \otimes \LL(V_{k,a_p,\chi})$.

Le théorème \ref{st} montre que $\Pibf(\dfont)^{\an}$ est dense dans $\Pibf(\dfont)$ et le résultat suivant, lui aussi montré dans \cite{CGL}, indique comment retrouver $\Pibf(\dfont)^{\an}$ en termes de la définition de $\Pibf(\dfont)$ donnée au (2) du théorème \ref{ccdf}.

\begin{theo}\label{locan}
Si $\dfont$ est un $(\phi,\Gamma)$-module étale sur $\calE$, absolument irréductible et de dimension $2$, alors $\Pibf(\dfont)^{\an}$ est l'image de $\dfont^\dag \boxtimes_\delta \Proj^1$ dans $\Pibf(\dfont) = \dfont \boxtimes_\delta \Proj^1 / (\dfont \boxtimes_\delta \Proj^1)^{\bor}$.
\end{theo}

Si $V$ est une représentation $E$-linéaire absolument irréductible de dimension $2$ de $\Gal(\Qpbar/\Qp)$ qui devient cristalline sur une extension abélienne de $\Qp$, alors Breuil avait conjecturé dans \cite{LB14} une description explicite de $\Pibf(V)^{\an}$ (comme deux induites paraboliques localement analytiques, amalgamées au-dessus d'une sous-représentation localement algébrique commune). Cette conjecture est démontrée dans \cite{LLA}.

\section{Applications}
\label{apps}

Nous donnons ici quelques applications de la correspondance et des propriétés de compatibilité qu'elle vérifie.

\subsection{Représentations triangulines}
\label{triapps}

Un sous-produit de la construction de la correspondance, et en particulier de la série principale unitaire, a été de dégager la notion de représentation $p$-adique trianguline de $\Gal(\Qpbar/\Qp)$. Comme on l'a dit au \S \ref{triang}, les représentations semi-stables sont \mbox{triangulines,} mais il en existe beaucoup d'autres. On a par exemple le résultat suivant de \cite{KOM} (pour les formes modulaires $p$-adiques et les représentations qui leur sont associées, voir le rapport \cite{EBS} d'Emerton).

\begin{theo}\label{ktr}
Les représentations associées aux formes modulaires paraboliques surconvergentes de pente finie sont triangulines.
\end{theo}

Tout comme la théorie de Hodge $p$-adique est un outil indispensable pour l'étude des formes modulaires, la théorie des représentations triangulines est au c{\oe}ur de l'étude des formes modulaires et automorphes $p$-adiques; c'est par exemple le thème de \cite{BCA} et des nombreux travaux qui s'en inspirent.

\subsection{Réduction des représentations cristallines}
\label{bug}

Soient $k \geq 2$ et $a_p \in \MM_E$ et $V_{k,a_p}$ la représentation définie au \S \ref{thpad}. Si on en choisit un $\OO_E$-réseau stable par $\Gal(\Qpbar/\Qp)$, que l'on réduit ce réseau modulo $\MM_E$ et qu'on semi-simplifie cette réduction, on obtient une représentation $k_E$-linéaire semi-simple $\overline{V}_{k,a_p}$ qui ne dépend pas du choix du réseau par le principe de Brauer-Nesbitt. La question se pose alors de déterminer $\overline{V}_{k,a_p}$. 

Si $k \leq p$, alors la réponse est donnée par la théorie de Fontaine-Laffaille (voir \cite{FL}) et on trouve que $\overline{V}_{k,a_p} = \ind(\omega_2^{k-1})$. En utilisant la théorie des modules de Wach (une spécialisation de la théorie des $(\phi,\Gamma)$-modules dans le cas cristallin), on peut calculer~$\overline{V}_{k,a_p}$ pour $k=p+1$ et pour $k \geq p+2$ si $\vp(a_p) \geq \lfloor (k-2) / (p+1) \rfloor$ (voir \cite{BLZ} pour ces calculs, et \cite{MV} pour des améliorations ponctuelles de la borne sur $\vp(a_p)$). Hors du disque $\vp(a_p) \geq \lfloor (k-2) / (p+1) \rfloor$, on n'a pas de formule générale et les calculs informatiques de \cite{SSB}, étendus ensuite par Buzzard, montrent que la situation se complique quand $k$ augmente. En particulier, $\overline{V}_{k,a_p}$ dépend de $a_p$ d'une manière de plus en plus compliquée. 

Le théorème \ref{compred} montre que $\overline{V}_{k,a_p}$ est déterminée par $\overline{\Pi}_{k,a_p}$, ce qui permet de calculer $\overline{V}_{k,a_p}$ dans les cas où on peut calculer $\overline{\Pi}_{k,a_p}$. C'est ce qu'a fait Breuil pour $k \leq 2p$ dans \cite{BR2} (et pour $k=2p+1$ dans un travail non publié) et qu'ont fait Buzzard et Gee dans \cite{BUG} pour $0 < \vp(a_p) < 1$. Le théorème suivant rassemble les résultats que l'on connaît pour l'instant (mars 2010).

\begin{theo}\label{allred}
La repr\'esentation $\overline{V}_{k,a_p}$ est connue dans les cas suivants
\begin{enumerate}
\item Si $2 \leq k \leq p+1$, alors $\overline{V}_{k,a_p} = \ind(\omega_2^{k-1})$.
\item Pour $k=p+2$ 
\begin{enumerate}
\item si $1 > \vp(a_p) > 0$, alors $\overline{V}_{k,a_p} = \ind(\omega_2^2)$
\item si $\vp(a_p) \geq 1$, et si $\lambda^2 - \overline{a_p/p} \cdot \lambda +1 = 0$, alors $\overline{V}_{k,a_p} = \omega\mu_{\lambda} \oplus \omega\mu_{\lambda^{-1}}$.
\end{enumerate}
\item Pour $2p \geq k \geq p+3$ 
\begin{enumerate}
\item si $1 > \vp(a_p) > 0$, alors $\overline{V}_{k,a_p} =   \ind(\omega_2^{k-p})$
\item si $\vp(a_p) = 1$, et si $\lambda=\overline{a_p/p} \cdot (k-1)$, alors 
$\overline{V}_{k,a_p} =\omega^{k-2} \mu_{\lambda} \oplus \omega\mu_{\lambda^{-1}}$
\item si $\vp(a_p) > 1$, alors $\overline{V}_{k,a_p} =   \ind(\omega_2^{k-1})$.
\end{enumerate}
\item Pour $k=2p+1$ (et $p \neq 2$) 
\begin{enumerate}
\item si $\vp(a_p^2+p) < 3/2$, alors $\overline{V}_{k,a_p} =   \ind(\omega_2^2)$
\item si $\vp(a_p^2+p) \geq 3/2$, alors $\overline{V}_{k,a_p} = \omega \mu_{\lambda} \oplus \omega\mu_{\lambda^{-1}}$ où $\lambda^2- \overline{(a_p^2+p)/(2 p a_p)} \cdot \lambda + 1 = 0$.
\end{enumerate}
\item Pour $k \geq 2p+2$, les r\'esultats ne sont que partiels
\begin{enumerate}
\item si $\vp(a_p)> \lfloor (k-2)/(p-1) \rfloor$, alors $\overline{V}_{k,a_p} =   \ind(\omega_2^{k-1})$ (qui est réductible si $p+1$ divise $k-1$)
\item si $0 < \vp(a_p) < 1$ et $t$ représente $k-1 \bmod{p-1}$ dans $\{1,\hdots,p-1\}$, alors 
\begin{enumerate}
\item $\overline{V}_{k,a_p} = \ind(\omega_2^t)$ si $p-1$ ne divise pas $k-3$
\item $\overline{V}_{k,a_p} \in \{ \ind(\omega_2^t)$,  $\omega \mu_{\lambda} \oplus \omega\mu_{\lambda^{-1}} \}$ pour un certain $\lambda$ si $p-1$ divise $k-3$.
\end{enumerate}
\end{enumerate}
\end{enumerate}
\end{theo}

On ne dispose pas pour l'instant de formule générale, même conjecturale, pour $\overline{V}_{k,a_p}$. Mentionnons tout de même la conjecture suivante de Buzzard (pour $p \neq 2$).

\begin{conj}
Si $k$ est pair et si $\overline{V}_{k,a_p}$ est réductible, alors $\vp(a_p)$ est un entier.
\end{conj}

Le problème analogue dans le cas semi-stable du calcul des $\overline{V}_{k,\Linv}$ se pose et le théorème \ref{compred} s'applique là aussi.

\subsection{Conjecture de Fontaine-Mazur}
\label{prfm}

Si $f = \sum_{n \geq 1} a_n q^n$ est une forme modulaire parabolique propre de poids $k$, 
de niveau~$N$ et de caractère $\varepsilon$, et si $E = \Qp(\{a_n\}_{n \geq 1})$, alors grâce à \cite{DMF} on sait qu'il existe une représentation $E$-linéaire $V_f$ de dimension $2$ de $\Gal(\Qbar/\QQ)$ telle que $V_f$ est non-ramifiée en tout $\ell \nmid pN$ et telle que, pour ces $\ell$, on a $\det(X-\mathrm{Fr}_\ell) = X^2-a_\ell X + \varepsilon(\ell)\ell^{k-1}$.

On sait par ailleurs que la restriction de $V_f$ à $\Gal(\Qpbar/\Qp)$ est potentiellement semi-stable et Fontaine et Mazur ont formulé dans \cite{FM} la conjecture suivante.

\begin{conj}\label{fm}
Si $V$ est une représentation $E$-linéaire irréductible de dimension $2$ de $\Gal(\Qbar/\QQ)$ qui est non-ramifiée en presque tout $\ell \neq p$, et dont la restriction à $\Gal(\Qpbar/\Qp)$ est potentiellement semi-stable à poids de Hodge-Tate distincts, alors il existe une forme modulaire parabolique propre telle que $V$ est la tordue de $V_f$ par un caractère.
\end{conj}

Notons $\overline{V}$ la réduction modulo $\MM_E$ de $V$.

\begin{theo}\label{emkis}
La conjecture de Fontaine-Mazur est vraie, si l'on suppose que $\overline{V}$ satisfait certaines hypothèses techniques.
\end{theo}

Ce théorème a été démontré indépendamment par Kisin (voir \cite{KFM}) et par Emerton (voir \cite{EFC}). Les {\og hypothèses techniques\fg} de Kisin sont les suivantes. 
\begin{enumerate}
\item $p \neq 2$ et $\overline{V}$ est impaire,
\item $\overline{V} \mid_{\Gal(\Qbar/\QQ(\zeta_p))}$ est absolument irréductible,
\item $\overline{V} \mid_{\Gal(\Qpbar/\Qp)}$ n'est pas de la forme $\smat{\omega \chi & * \\ 0 & \chi}$.
\end{enumerate}

Les {\og hypothèses techniques\fg} d'Emerton sont (1) et (2) et
\begin{itemize}
\item[3'.] $\overline{V} \mid_{\Gal(\Qpbar/\Qp)}$ n'est pas de la forme $\smat{\chi & * \\ 0 & \omega \chi}$ ni de la forme $\smat{\chi & * \\ 0 & \chi}$.
\end{itemize}

La méthode d'Emerton donne alors un résultat supplémentaire : si on suppose que $V \mid_{\Gal(\Qpbar/\Qp)}$ est trianguline (au lieu de potentiellement semi-stable), alors $V$ provient d'une forme modulaire parabolique surconvergente de pente finie.

L'outil le plus puissant dont on dispose pour l'instant afin de démontrer la modularité de certaines représentations galoisiennes est l'étude de leurs espaces de déformations. C'est cette méthode qui a permis à Wiles de démontrer dans \cite{WST} la modularité des courbes elliptiques semi-stables. 

Dans leur article \cite{BM1}, Breuil et Mézard ont proposé une conjecture reliant certains anneaux paramétrant les déformations potentiellement semi-stables d'une représentation $\overline{V}$, et certaines représentations de $\GL_2(\Zp)$. Plus précisément, pour un certain type de paramètres $(k,\tau,\overline{V})$ de déformations, ils définissent une multiplicité galoisienne $\mu_{\Gal}(k,\tau,\overline{V})$ et une multiplicité automorphe $\mu_{\Aut}(k,\tau,\overline{V})$. La multiplicité galoisienne mesure la {\og taille\fg} de l'anneau des déformations de type $(k,\tau,\overline{V})$, tandis que la multiplicité automorphe dépend de $\overline{V}$ et de la réduction modulo $p$ de certaines représentations de $\GL_2(\Zp)$ associées à $k$ et $\tau$ par la correspondance locale de Langlands {\og classique\fg}. La conjecture de Breuil-Mézard est alors que $\mu_{\Gal}(k,\tau,\overline{V}) = \mu_{\Aut}(k,\tau,\overline{V})$.

Le foncteur de Colmez $V \mapsto \Pibf(V)$ étant défini de manière assez naturelle, il s'étend aux familles et définit par suite un foncteur de l'espace des déformations de $\overline{V}$ vers l'espace des déformations de $\overline{\Pibf}(V)$ (c'est d'ailleurs cet argument qui est utilisé par Kisin dans \cite{KAC}). Cette construction, ainsi que le théorème \ref{localg} concernant les vecteurs localement algébriques, permettent à Kisin de faire le lien entre multiplicité galoisienne et multiplicité automorphe, et par suite de démontrer la conjecture de Breuil-Mézard. Ceci lui donne des renseignements précis sur les anneaux de déformations de $\overline{V}$ qui lui permettent alors d'appliquer les techniques de modularité et de démontrer la conjecture de Fontaine-Mazur.

Remarquons pour terminer que la conjecture de Breuil-Mézard était inspirée des calculs de \cite{BCDT} dont le principal résultat était la modularité de toutes les courbes elliptiques définies sur $\QQ$ et c'est la rédaction de \cite{BM1} qui a contribué à donner à Breuil l'idée de la construction de la correspondance de Langlands $p$-adique. La boucle est donc bouclée, puisque cette correspondance de Langlands $p$-adique pour $\GL_2(\Qp)$ permet à présent de démontrer une vaste généralisation du résultat de \cite{BCDT}.

\end{document}